\documentclass[12pt]{article}
\usepackage{amsmath,amssymb,latexsym}
\usepackage{amsthm}
\usepackage[dvips]{graphics}
\newtheorem{thrm}{Theorem}[section]
\newtheorem{lemm}{Lemma}[section]
\newtheorem{dfnt}{Definition}[section]
\newtheorem{prot}{Proposition}[section]

\newtheorem{corl}{Corollary}[section]
\newtheorem{remk}{Remark}[section]
\newtheorem{exa}{Example}[section]
\newtheorem{nte}{Note}[section]
\numberwithin{equation}{section}

\usepackage{lineno}
\begin{document}
\title{On linearity of pan-integral and pan-integrable \\
                           functions space}
\author{
Yao Ouyang$\sp{a}$\thanks{Corresponding author. E-mail:oyy@hutc.zj.cn(Yao Ouyang); lijun@cuc.edu.cn(Jun Li); mesiar@math.sk(Radko Mesiar)}\qquad Jun Li$\sp{b}$\qquad Radko Mesiar$\sp{c,d}$ \\
\small{\it $\sp{a}$Faculty of Science, Huzhou Teachers College, Huzhou, Zhejiang 313000, China}\\
\small\it$\sp{b}$Department of Applied Mathematics, College of Science, Communication University of China,\\
      \small\it{Beijing 100024, People's Republic of China }\\
\small{\it$\sp{c}$ Slovak University of Technology, Faculty of Civil Engineering,}\\
       \small{\it Radlinsk\' eho 11, 811 05 Bratislava, Slovakia}\\
\small{\it $\sp{d}$ UTIA CAS, Pod Vod\'arenskou v\v e\v z\'i 4, 182 08 Prague, Czech Republic}
}
\date{}
\maketitle

\begin{abstract}
$L\sp{p}$ space is a crucial aspect of classical measure theory. For nonadditive measure, it is known that $L\sp{p}$ space theory holds for the Choquet integral whenever the monotone measure $\mu$ is submodular and continuous from below. The main purpose of this paper is to generalize $L\sp{p}$ space theory to $+,\cdot$-based pan-integral. Let $(X, {\cal A}, \mu)$ be a monotone measure space. We prove that the $+,\cdot$-based pan-integral is additive with respect to integrands if $\mu$ is subadditive. Then we introduce the pan-integral for real-valued functions(not necessarily nonnegative), and prove that this integral possesses linearity if $\mu$ is subadditive. By using the linearity of pan-integral, we finally show that all of the pan-integrable functions form a Banach space. Since the $+,\cdot$-based pan-integral coincides with the concave integral for subadditive measure, the results obtained in this paper remain valid for the concave integral. Noticing that an outer measure is subadditive, we can define a Lebesgue-like integral(possesses linearity) from an outer measure, and the $L\sp{p}$ theory holds for this integral.
\end{abstract}
\section{Introduction}
Let $(X, {\cal A}, m)$ be a (classical) measure space and $E$ be a measurable subset of $X$. By ${\cal F}(E)$ we denote the set of all measurable functions defined on $E$. Then for any $p\geq 1$, the set
\[
\left\{f\in {\cal F}(E)\colon\left(\int|f|\sp{p} dm\right)\sp{\frac{1}{p}}<\infty\right\}
\]
forms a complete and separable metric space. This is the well-known $L\sp{p}$ space, which is an important aspect of the Lebesgue integral theory. The linearity of the Lebesgue integral plays a crucial role in the $L\sp{p}$ space theory.

The Choquet integral, as a prominent nonlinear integral, was introduced by Choquet \cite{Cho53}. The $L\sp{p}$ space theory is not true for the Choquet integral in general case, since this integral lacks additivity(in fact, the Choquet integral is a level set-based integral, which ensures its comonotonic additivity- a property weaker than additivity). Denneberg \cite{Den94} showed a subadditivity theorem for the Choquet integral under the assumption that the monotone measure is submodular. Then the $L\sp{p}$ space theory was generalized to Choquet integral (see \cite{Den94}). The recent paper by Shirali \cite{Shi08} presented a similar result.

The pan-integral, which was introduced by Yang \cite{Yan85}, is another well-known nonlinear integral. Unlike the Choquet integral, the pan-integral is a partition-based integral. As a result, the pan-integral even lacks comonotonic additivity. From this point of view, it requires to exam whether or not the $L\sp{p}$ space theory holds for the pan-integral.

The paper is structured as follows. After this introduction, we recall some basic facts about the monotone measure, the Choquet and pan-integrals in Section 2. Sections 3 and 4 discuss the properties of the pan-integral for nonnegative functions, the main point in which is the additivity (w.r.t. the integrands)of a pan-integral if the monotone measure is subadditive. Section 5 defines the pan-integral for general real-valued functions. This integral is introduced in the same way as the symmetric Choquet integral, and it possesses linearity if the monotone measure is subadditive. Section 6 generalizes the $L\sp{p}$ space theory to the pan-integral and Section 7 concludes.

\section{Preliminaries}
Let $X$ be a nonempty set and ${\mathcal A}$ a $\sigma$-algebra of subsets of $X$.
${\cal F}$ denotes the class of all real-valued measurable functions
on $(X, {\mathcal A})$, and ${\cal F}\sp{+}$ is the collection of those nonnegative functions
which are taken from ${\cal F}$. Unless stated otherwise all the
subsets mentioned are supposed to belong to ${\mathcal A}$ and all functions appeared in this paper are supposed to belong to ${\cal F}$.

By a finite partition of $(X, {\mathcal A})$ means that it is a finite disjoint system of sets $\{A\sb{i}\}\sb{i=1}\sp{n}\subset{\mathcal A}$ such that
$A\sb{i}\cap A\sb{j}=\emptyset$ for $i\neq j$ and $\cup\sb{i=1}\sp{n} A\sb{i}=X$.
Let $\hat{\cal{P}}$ be the set of all finite partitions of
$(X, {\mathcal A})$.

We assume that $\mu$ is a monotone measure on $(X, {\cal A})$,
 i.e.,
 $\mu: {\mathcal{A}} \to [0, +\infty]$ is an extended real valued set function
satisfying the following conditions:

{\rm (1)} \ \ $\mu(\emptyset)=0$ and $\mu(X) > 0$;

{\rm (2)} \ \ $\mu(A) \leq \mu(B)$ whenever $A \subset B$ and $A, B
\in {\mathcal{A}}$.

When $\mu$ is a monotone measure, the triple $(X, {\mathcal A}, \mu)$ is called
a monotone measure space (\cite{Pap95,WanKli09}).
In the literature, the monotone measure constraint by $\mu(X)=1$ (sometimes also need the continuity(see below) of $\mu$) is also known as a capacity or
fuzzy measure, or nonadditive probability, etc. (see \cite{Den94,Leh09,Ogi94,Vic08,WanWanKli96}).

 Recall that a set function $\mu$ is said to be

 (i) \ \emph{continuous from below}, if for any $\{E\sb{n}\}\sb{n=1}\sp{\infty}\subset {\mathcal A}$, $E\sb{1}\subset E\sb{2}\subset\dots$ implies $\mu(\cup\sb{n=1}\sp{\infty} E\sb{n})=\lim\sb{n\to\infty}\mu(E\sb{n})$;

 (ii) \ \emph{continuous from above}, if for any $\{E\sb{n}\}\sb{n=1}\sp{\infty}\subset {\mathcal A}$, $E\sb{1}\supset E\sb{2}\supset\dots$ and $\mu(E\sb{1})<\infty$ imply $\mu(\cap\sb{n=1}\sp{\infty} E\sb{n})=\lim\sb{n\to\infty}\mu(E\sb{n})$;

 (iii) \emph{continuous} if it is continuous both from below and from above;

 (iv)  \emph{submodular}, if $\mu(A \cup B)+ \mu(A\cap B)\leq  \mu(A) + \mu(B)$ for any $A, B \in {\mathcal A}$;

 (v) \emph{supermodular}, if
 $\mu(A \cup B)+ \mu(A\cap B)\geq  \mu(A) + \mu(B)$ for any $A, B \in {\mathcal A}$;

 (vi) \emph{subadditive}, if $\mu(A \cup B)\leq  \mu(A) + \mu(B)$ for any $A, B \in {\mathcal A}$;

 (vii) \emph{null-additive}, if for any $A, B \in {\mathcal A}$, $\mu(B)=0$ implies $\mu(A\cup B)=\mu(A)$.

 Obviously, the submodularity of $\mu$ implies the subadditivity and both imply the null-additivity, but not vice versa.

In \cite{WanKli09} (see also \cite{Yan85}) the concept of pan-integral was introduced, which involves two binary operations, the pan-addition $\oplus$ and
pan-multiplication $\otimes$ of real numbers (see \cite{Pap95,SugMur87,WanKli09,WanWanKli96}).
In this paper we only consider the pan-integrals based on
the standard addition $+$ and multiplication $\cdot$. We recall the following definition.


%
\begin{dfnt}
   \label{Def-31}
{\em

Let $(X, {\mathcal A}, \mu)$ be a monotone measure space, $A\in \mathcal A$ and $f \in {\cal F}\sp{+}$.
The pan-integral of $f$ on $A$ with respect to $\mu$, is defined by
    $$
   \int_A\sp{pan} f d\mu
   =\sup\limits_{\mathcal{E} \in {\hat {\cal P}}} \left\{\sum\limits_{E \in \mathcal{E}}
    \big[(\inf\limits_{x \in A\cap E} f(x))\cdot \mu(A\cap E)\big]\right\}.
    \eqno{(2.1)}
    $$
When $A=X$,  $\int_X\sp{pan} f d\mu $ will be denoted by $\int\sp{pan} fd\mu$ for short.
    }
\end{dfnt}

The above pan-integral of $f$ can be expressed as the following form:
\begin{eqnarray*}
\int\sp{pan} f d\mu
       &=& \sup\bigg\{ \sum\sb{i=1}\sp{n}\lambda\sb{i}\mu(A\sb{i}): \sum\sb{i=1}\sp{n}\lambda\sb{i}\chi\sb{A\sb{i}}\leq f, \\[-2mm]
       & & \hspace{12mm} \{A\sb{i}\}\sb{i=1}\sp{n}\subset {\mathcal A}
                  \mbox{ is a partition of } X, \lambda\sb{i}\geq 0,
                        n\in \mathbb{N} \bigg\},
\end{eqnarray*}
where $\chi\sb{A_{i}}$ is the characteristic function of $A_{i}$.

In order to simplify our discussions, we introduce the concept positive set of a nonnegative function.
\begin{dfnt} {\em
Let $f\in {\cal F}\sp{+}$. The set $\{x\in X: f(x)> 0\}$ (briefly $\{f> 0\}$)
is called {\it positive set of f}. 
For $f, g\in {\cal F}\sp{+}$, if $\{f> 0\}\cap\{g> 0\}=\emptyset$, then we say that
the positive sets of $f$ and $g$ are disjoint.
}
\end{dfnt}

In the following we present some basic properties of pan-integrals, some of them can be found in \cite{WanKli09}.
\begin{prot}\label{prot2-3}
Let $(X, {\cal A}, \mu)$ be a monotone measure space and $f, g\in {\cal F}\sp{+}$. Then we have the following

(i) if $\mu(A)=0$ then $\int\sb{A}\sp{pan} f d\mu=0$;

(ii) if $f=0$ ~$\mu$-a.e. on $A$ , i.e.,  $\mu(\{ x\in A: f(x)\neq 0 \})=0$,
then $\int\sb{A}\sp{pan} f d\mu=0$. Forthmore, if $\mu$ is continuous from below, then $\int\sb{A}\sp{pan} f d\mu=0$ if and only if $f=0$ ~$\mu$-a.e. on $A$;

(iii) if $f\leq g$ on $A$, then $\int\sb{A}\sp{pan} f d\mu\leq \int\sb{A}\sp{pan} g d\mu$;

(iv)  for  $A\in {\mathcal A}$, $\int\sb{A}\sp{pan} f d\mu=\int\sp{pan} f\cdot\chi\sb{A}d\mu$.

(v) $\int\sp{pan} f d\mu=\int\sb{\{f>0\}}\sp{pan} f d\mu$.

(vi) If $\mu$ is null-additive and $f=g$~ $\mu$-a.e. on $A$, then  $\int\sb{A}\sp{pan} f d\mu= \int\sb{A}\sp{pan} g d\mu$.
\end{prot}

{\bf Proof.} We only prove (vi). Let $A\sb{1}=\{x\in A, f(x)\neq g(x)\}$, then $\mu(A\sb{1})=0$ and for any subset $B$ of $A$ we have $\mu(B)=\mu(B\setminus A\sb{1})$. Thus, for any $\sum\sb{i=1}\sp{n} \lambda\sb{i}\chi\sb{B\sb{i}}\leq f(x)\cdot\chi\sb{A}$, $\sum\sb{i=1}\sp{n} \lambda\sb{i}\chi\sb{B\sb{i}\setminus A\sb{1}}\leq g(x)\cdot\chi\sb{A}$. Moreover,
\[
\int\sb{A}\sp{pan} gd\mu\geq\sum\sb{i=1}\sp{n} \lambda\sb{i}\mu(B\sb{i}\setminus A\sb{1})=\sum\sb{i=1}\sp{n} \lambda\sb{i}\mu(B\sb{i}),
\]
which implies that $\int\sb{A}\sp{pan} gd\mu\geq \int\sb{A}\sp{pan} fd\mu$. In a similar way $\int\sb{A}\sp{pan} fd\mu\geq \int\sb{A}\sp{pan} gd\mu$, and thus (vi) holds.$\Box$
\begin{nte}\label{nte2-4}{\em
(1) When $\mu$ is continuous from below, from the above proposition (ii), we know that for any measurable functions $f, g\in {\cal F}$ and any real number $p>0$, $\int\sb{A}\sp{pan} |f-g|\sp{p} d\mu=0$ if and only if $f=g$ ~$\mu$-a.e. on $A$.

(2) Due to the above proposition (iv), we will not distinguish the partition of $X$ and the partition of $\{f>0\}$. Sometimes we may just say that $\{A\sb{i}\}\sb{i=1}\sp{n}$ is a partition.

(3) By (vi), if $\mu$ is null-additive then, from the pan-integral point of view, we do not distinguish $f$ and $g$ whenever $f=g$~ $\mu$-a.e. on $X$.
 }
\end{nte}
\section{The additivity of the pan-integral for nonnegative measurable functions}
The main task of this section is to prove the following result.

\begin{thrm}\label{th3-1}
Let $(X, {\cal A
}, \mu)$ be a monotone measure space. If $\mu$ is subadditive, then the pan-integral is additive w.r.t. the integrands, i.e., for any $f, g\in {\cal F\sp{+}}$,
\begin{equation}\label{eq3-1}
\int\sp{pan}(f+g) d\mu=\int\sp{pan} f d\mu+ \int\sp{pan} g d\mu
\end{equation}
\end{thrm}

{\bf Proof.} Here we only consider the case that $f+g$ is pan-integrable, that is, $\int\sp{pan} (f+g) d\mu$ is finite. The case of $\int\sp{pan} (f+g) d\mu=\infty$ will be considered in the next section.

For an arbitrary but fixed $\varepsilon>0$, there exist a partition $\{A\sb{i}\}\sb{i=1}\sp{n}$ and a sequence of nonnegative numbers $\{\lambda\sb{i}\}\sb{i=1}\sp{n}$ such that both $\sum\sb{i=1}\sp{n}\lambda\sb{i}\chi\sb{A\sb{i}}\leq f+g$ and
\[
\int\sp{pan}\left(f+g\right) d\mu<\sum\sb{i=1}\sp{n}\lambda\sb{i}\mu(A\sb{i})+\frac{\varepsilon}{2}
\]
hold. Without loss of generality, we can assume that $\lambda\sb{i}>0$ for all $i$ (otherwise, the summand $\lambda\sb{i}\mu(A\sb{i})$ can be omitted). Since $\mu(A\sb{1})<\infty$ and $\lambda\sb{1}>0$, there exists a positive number $l\sb{1}<\lambda\sb{1}$ such that $(\lambda\sb{1}-l\sb{1})\mu(A\sb{1})<\frac{\varepsilon}{2\sp{2}}$. Obviously, $l\sb{1}\chi\sb{A\sb{1}}<\lambda\sb{1}\chi\sb{A\sb{1}}\leq(f+g)\cdot\chi\sb{A\sb{1}}$. Let $\delta=\lambda\sb{1}-l\sb{1}$ and divide the interval $[0, l\sb{1}]$ as
\[
0=r\sb{1}\sp{(0)}<r\sb{1}\sp{(1)}<\dots<r\sb{1}\sp{(m\sb{1})}=l\sb{1}
\]
such that $\max\sb{1\leq i\leq m\sb{1}}(r\sb{1}\sp{(i)}-r\sb{1}\sp{(i-1)})<\delta$. Denote
\[
A\sb{1}\sp{(k)}=\left\{x\in A\sb{1}|r\sb{1}\sp{(k-1)}\leq f(x)<r\sb{1}\sp{(k)}\right\}, k=1,\dots,m\sb{1},
\]
and
\[
A\sb{1}\sp{(m\sb{1}+1)}=\left\{x\in A\sb{1}| f(x)\geq l\sb{1}\right\}.
\]
Then $A\sb{1}\sp{(i)}\cap A\sb{1}\sp{(j)}=\emptyset (i\neq j)$ and $\cup\sb{i=1}\sp{m\sb{1}+1}A\sb{1}\sp{(i)}=A\sb{1}$. Moreover, we conclude that for each $x\in A\sb{1}\sp{(k)}, k=1,2,\dots,m\sb{1}$, $g(x)\geq l\sb{1}-r\sb{1}\sp{(k-1)}$. In fact, if $g(x)< l\sb{1}-r\sb{1}\sp{(k-1)}$, then
\[
f(x)+g(x)<r\sb{1}\sp{(k)}+(l\sb{1}-r\sb{1}\sp{(k-1)})=l\sb{1}+(r\sb{1}\sp{(k)}-r\sb{1}\sp{(k-1)})<l\sb{1}+\delta=\lambda\sb{1},
\]
which contradicts with the fact that $\lambda\sb{1}\chi\sb{A\sb{1}}\leq (f+g)\cdot\chi\sb{A\sb{1}}$. For $x\in A\sb{1}\sp{(m\sb{1}+1)}$, $g(x)\geq 0=l\sb{1}-r\sb{1}\sp{(m\sb{1})}$ also holds. Denote
\[
t\sb{1}\sp{(k)}=l\sb{1}-r\sb{1}\sp{(k-1)}, k=1,\dots,m\sb{1}+1.
\]
Then
\[
\sum\sb{k=1}\sp{m\sb{1}+1}r\sb{1}\sp{(k-1)}\chi\sb{A\sb{1}\sp{(k)}}\leq f\cdot\chi\sb{A\sb{1}}
\]
and
\[
\sum\sb{k=1}\sp{m\sb{1}+1}t\sb{1}\sp{(k)}\chi\sb{A\sb{1}\sp{(k)}}\leq g\cdot\chi\sb{A\sb{1}}.
\]
Observe that the subadditivity of $\mu$ implies
\[
\mu(A\sb{1})=\mu\left(\bigcup\sb{k=1}\sp{m\sb{1}+1}A\sb{1}\sp{(k)}\right)\leq\sum\sb{k=1}\sp{m\sb{1}+1}\mu\left(A\sb{1}\sp{(k)}\right).
\]
Thus
\begin{eqnarray*}
&&\sum\sb{k=1}\sp{m\sb{1}+1}r\sb{1}\sp{(k-1)}\mu\left(A\sb{1}\sp{(k)}\right)+\sum\sb{k=1}\sp{m\sb{1}+1}t\sb{1}\sp{(k)}\mu\left(A\sb{1}\sp{(k)}\right)\\
&=&\sum\sb{k=1}\sp{m\sb{1}+1}\left(r\sb{1}\sp{(k-1)}+t\sb{1}\sp{(k)}\right)\mu\left(A\sb{1}\sp{(k)}\right)\\
&=&\sum\sb{k=1}\sp{m\sb{1}+1}l\sb{1}\mu\left(A\sb{1}\sp{(k)}\right)\geq l\sb{1}\mu(A\sb{1})\\
&>& \lambda\sb{1}\mu(A\sb{1})-\frac{\varepsilon}{2\sp{2}}.
\end{eqnarray*}

For each $i\leq n$, let $l\sb{i}\in (0, \lambda\sb{i})$ be such that $(\lambda\sb{i}-l\sb{i})\mu(A\sb{i})<\frac{\varepsilon}{2\sp{i+1}}$. By using the technique used above, we can prove that for each $j\leq n$ there exist a partition $\{A\sb{j}\sp{(k)}\}\sb{k=1}\sp{m\sb{j}+1}$ of $A\sb{j}$ and two sequences of nonnegative number $\{r\sb{j}\sp{(k-1)}\}\sb{k=1}\sp{m\sb{j}+1}$ and $\{t\sb{j}\sp{(k)}=l\sb{j}-r\sb{j}\sp{(k-1)}\}\sb{k=1}\sp{m\sb{j}+1}$ with
\[
\sum\sb{k=1}\sp{m\sb{j}+1}r\sb{j}\sp{(k-1)}\chi\sb{A\sb{j}\sp{(k)}}\leq f\cdot\chi\sb{A\sb{j}}
\]
and
\[
\sum\sb{k=1}\sp{m\sb{j}+1}t\sb{j}\sp{(k)}\chi\sb{A\sb{j}\sp{(k)}}\leq g\cdot\chi\sb{A\sb{j}},
\]
such that
\[
\sum\sb{k=1}\sp{m\sb{j}+1}r\sb{j}\sp{(k-1)}\mu\left(A\sb{j}\sp{(k)}\right)+\sum\sb{k=1}\sp{m\sb{j}+1}t\sb{j}\sp{(k)}\mu\left(A\sb{j}\sp{(k)}\right)
>\lambda\sb{j}\mu(A\sb{j})-\frac{\varepsilon}{2\sp{j+1}}.
\]
Thus, we have proven that
\[
\sum\sb{j=1}\sp{n}\sum\sb{k=1}\sp{m\sb{j}+1}r\sb{j}\sp{(k-1)}\chi\sb{A\sb{j}\sp{(k)}}\leq\sum\sb{j=1}\sp{n} f\cdot\chi\sb{A\sb{j}}\leq f,
\]
\[
\sum\sb{j=1}\sp{n}\sum\sb{k=1}\sp{m\sb{j}+1}t\sb{j}\sp{(k)}\chi\sb{A\sb{j}\sp{(k)}}\leq\sum\sb{j=1}\sp{n} g\cdot\chi\sb{A\sb{j}}\leq g,
\]
and
\begin{eqnarray*}
&&\int\sp{pan}f d\mu+\int\sp{pan}g d\mu\\
&\geq&\sum\sb{j=1}\sp{n}\sum\sb{k=1}\sp{m\sb{j}+1}r\sb{j}\sp{(k-1)}\mu\left(A\sb{j}\sp{(k)}\right)
+\sum\sb{j=1}\sp{n}\sum\sb{k=1}\sp{m\sb{j}+1}t\sb{j}\sp{(k)}\mu\left(A\sb{j}\sp{(k)}\right)\\
&>& \sum\sb{j=1}\sp{n}\left(\lambda\sb{j}\mu(A\sb{j})-\frac{\varepsilon}{2\sp{j+1}}\right)\\
&>& \sum\sb{j=1}\sp{n}\lambda\sb{j}\mu(A\sb{j})-\frac{\varepsilon}{2}\\
&>& \int\sp{pan}(f+g) d\mu-\varepsilon.
\end{eqnarray*}
Letting $\varepsilon\to 0$, we get $\int\sp{pan}fd\mu+\int\sp{pan}gd\mu\geq \int\sp{pan}(f+g)d\mu$ as desired.

On the other hand, for any $\sum\sb{i=1}\sp{n} \alpha\sb{i}\chi\sb{A\sb{i}}\leq f$ and $\sum\sb{j=1}\sp{m} \beta\sb{j}\chi\sb{B\sb{j}}\leq g$, where $\alpha\sb{i}, \beta\sb{j}\geq 0$, $\{A\sb{i}\}$ and $\{B\sb{j}\}$ are two partitions of $X$, put $C\sb{ij}=A\sb{i}\cap B\sb{j}$ ($C\sb{ij}=\emptyset$ may hold for some $i, j$) and $\gamma\sb{ij}=\alpha\sb{i}+\beta\sb{j}$, $1\leq i\leq n, 1\leq j\leq m$. Then
\begin{eqnarray*}
\sum\sb{i=1}\sp{n}\sum\sb{j=1}\sp{m} \gamma\sb{ij}\chi\sb{C\sb{ij}}
&=&
\sum\sb{i=1}\sp{n}\sum\sb{j=1}\sp{m}\alpha\sb{i}\chi\sb{C\sb{ij}}+\sum\sb{j=1}\sp{m}\sum\sb{i=1}\sp{n}\beta\sb{j}\chi\sb{C\sb{ij}}\\
&=&\sum\sb{i=1}\sp{n}\alpha\sb{i}\chi\sb{A\sb{i}}+\sum\sb{j=1}\sp{m}\beta\sb{j}\chi\sb{B\sb{j}}\leq f+g
\end{eqnarray*}
and
\begin{eqnarray*}
\int\sp{pan} (f+g) d\mu
&\geq&
\sum\sb{i=1}\sp{n}\sum\sb{j=1}\sp{m} \gamma\sb{ij}\mu(C\sb{ij})\\
&=&
\sum\sb{i=1}\sp{n}\sum\sb{j=1}\sp{m} \alpha\sb{i}\mu(C\sb{ij})+ \sum\sb{j=1}\sp{m}\sum\sb{i=1}\sp{n} \beta\sb{j}\mu(C\sb{ij})\\
&\geq&
\sum\sb{i=1}\sp{n} \alpha\sb{i}\mu\Big(\bigcup\sb{j=1}\sp{m} C\sb{ij}\Big)+ \sum\sb{j=1}\sp{m} \beta\sb{j}\mu\Big(\bigcup\sb{i=1}\sp{n} C\sb{ij}\Big)\\
&=&
\sum\sb{i=1}\sp{n} \alpha\sb{i}\mu(A\sb{i})+ \sum\sb{j=1}\sp{m} \beta\sb{j}\mu(B\sb{j}),
\end{eqnarray*}
from which we can get that
\[
\int\sp{pan}(f+g) d\mu\geq \int\sp{pan}f d\mu+ \int\sp{pan}g d\mu.
\]

Hence it holds that $\int\sp{pan}(f+g) d\mu= \int\sp{pan}f d\mu+ \int\sp{pan}g d\mu$.$\Box$

\begin{nte}
Theorem \ref{th3-1} tells us that if $\mu$ is subadditive then all nonnegative pan-integrable functions form a convex cone.
\end{nte}

\begin{corl}
Let $(X, {\cal A}, \mu)$ be a monotone measure space. If $\mu$ is subadditive, then for any $f\in{\cal F\sp{+}}$ and any $A, B\in{\cal A}$ with $A\cap B=\emptyset$, we have
\[
\int\sp{pan}\sb{A\cup B}fd\mu=\int\sp{pan}\sb{A}fd\mu+ \int\sp{pan}\sb{B}fd\mu.\Box
\]
\end{corl}

{\bf Proof.}~Denote $g=f\cdot\chi\sb{A}$ and $h=f\cdot\chi\sb{B}$. Then, by (iv) of Proposition \ref{prot2-3} and Theorem \ref{th3-1}, we have
\begin{eqnarray*}
\int\sp{pan}\sb{A\cup B}fd\mu &=& \int\sp{pan}f\cdot\chi\sb{A\cup B}d\mu\\
&=& \int\sp{pan} (g+h) d\mu\\
&=&\int\sp{pan} g d\mu+ \int\sp{pan} h d\mu\\
&=& \int\sp{pan}\sb{A}fd\mu+ \int\sp{pan}\sb{B}fd\mu.
\end{eqnarray*}

When $X$ is finite, Theorem \ref{th3-1} has obvious meaning. Let $(X, {\cal F}, \mu)$ be a monotone measure space with $|X|<\infty$ and $A\sb{1},\dots,A\sb{k}$ be minimal atoms. Recall that a minimal atom is a measurable set $A$ such that $\mu(A)>0$ and for each measurable proper subset $B$ of $A$ we have $\mu(B)=0$. When $X$ is a finite set, each set of positive measure contains at least one minimal atom \cite{OuyLiMes15}. If we further assume that $\mu$ is subadditive then for each minimal atom $A$, there are no nonempty proper subsets $B$ of $A$ such that $B\in{\cal F}$. In fact, if there are some nonempty proper subset $B\in{\cal F}$, then by the definition of minimal atom, $\mu(B)=0, \mu(A\setminus B)=0$ and the subadditivity of $\mu$ implies $\mu(A)=0$, a contradiction. Observing this fact, to ensure its measurability, a function $f$ must take constant value on each minimal atom. Thus, for any measurable function $f\colon X\to [0, \infty)$ it holds
\[
\int\sp{pan}fd\mu=\sum\sb{i=1}\sp{k} c\sb{i}\mu(A\sb{i}),
\]
where $c\sb{i}$ is the constant value $f$ takes on the minimal atom $A\sb{i}$. From this fact it is immediate that the pan-integral is additive.

\section{Further discussion of Theorem \ref{th3-1}}

In this section, we reconsider Theorem \ref{th3-1} for the case that $f+g$ is not pan-integrable. We need two lemmas.

\begin{lemm}\label{lemm4-1}
Let $(X, {\cal A}, \mu)$ be a monotone measure space, and $f, g\in {\cal F\sp{+}}$.
If $\{f>0\}\cap~\{g>0\}=\emptyset$, then
\begin{equation}\label{eq4-1}
\int\sp{pan}(f+g)d\mu\geq \int\sp{pan}fd\mu+\int\sp{pan}gd\mu.
\end{equation}
\end{lemm}

{\bf Proof.~} If one of the two integrals on the right-hand side of Ineq. (\ref{eq4-1}) is infinite then, by the monotonicity of the pan-integral, $\int\sp{pan}(f+g)d\mu$ is also equal to infinity, which implies the validity of (\ref{eq4-1}).

So, without loss of generality, we can suppose that both $\int\sp{pan}fd\mu$ and $\int\sp{pan}gd\mu$ are finite. For any arbitrary but fixed $\varepsilon>0$, there are a partition $\{A\sb{i}\}\sb{i=1}\sp{k}$ of $\{f>0\}$, a partition $\{B\sb{j}\}\sb{j=1}\sp{m}$ of $\{g>0\}$, and two sequences of positive numbers $\{\lambda\sb{i}\}\sb{i=1}\sp{k}$ and $\{l\sb{j}\}\sb{j=1}\sp{m}$ such that $\sum\sb{i=1}\sp{k}\lambda\sb{i}\chi\sb{A\sb{i}}\leq f$,  $\sum\sb{j=1}\sp{m}l\sb{j}\chi\sb{B\sb{j}}\leq g$ and both the following two inequalities hold
\[
\int\sp{pan} fd\mu<\sum\sb{i=1}\sp{k}\lambda\sb{i}\mu(A\sb{i})+\frac{\varepsilon}{2},\qquad \int\sp{pan} gd\mu<\sum\sb{j=1}\sp{m}l\sb{j}\mu(B\sb{j})+\frac{\varepsilon}{2}.
\]
By the fact of $\{f>0\}\cap \{g>0\}=\emptyset$, we know that $\{A\sb{i}\}\sb{i=1}\sp{k}\cup\{B\sb{j}\}\sb{j=1}\sp{m}$ is a partition of $\{f+g>0\}$. Moreover, we have that $\sum\sb{i=1}\sp{k}\lambda\sb{i}\chi\sb{A\sb{i}}+\sum\sb{j=1}\sp{m}l\sb{j}\chi\sb{B\sb{j}}\leq f+g$, and
\begin{eqnarray*}
\int\sp{pan}(f+g)d\mu
&\geq&
\sum\sb{i=1}\sp{k}\lambda\sb{i}\mu(A\sb{i})+ \sum\sb{j=1}\sp{m}l\sb{j}\mu(B\sb{j})\\
&\geq& \int\sp{pan} fd\mu+ \int\sp{pan} gd\mu-\varepsilon.
\end{eqnarray*}
Letting $\varepsilon\to 0$, we get the result as desired.$\Box$

\begin{nte}
  The property described in Lemma \ref{lemm4-1} can be called the positive sets disjointness superadditivity of the pan-integral. In \cite{OuyLilnai}, we also call this property the support disjointness superadditivity. The following example shows that Ineq. (\ref{eq4-1}) can be violated if $\{f>0\}\cap~\{g>0\}\neq\emptyset$.
\end{nte}

\begin{exa}

Let $X=\{1,2,3,\dots\}, {\cal F}=2\sp{X}$ and the monotone measure $\mu\colon {\cal F}\to [0, \infty]$ be defined as $\mu(A)=\frac{1}{|X\setminus A|+1}$, where $|A|$ stands for the cardinality of $A$. Suppose that $f, g\colon X\to [0, \infty]$ are defined as $f(x)=x-1$ and $g(x)=1$ for all $x\in X$. Then, we can show that
\[
\int\sp{pan} (f+g)d\mu<\int\sp{pan} fd\mu+\int\sp{pan} gd\mu.
\]
Observe first that for any subset $A$, $\mu(A)>0$ if and only if $X\setminus A$ is a finite set. So, for any partition $\{A\sb{i}\}\sb{i=1}\sp{k}$ of $X$, there are at most one set $A\sb{i}$ with positive measure. For simplicity, we suppose that $\mu(A\sb{1})>0$ and $t=\min\{x|x\in A\sb{1}\}$. Then $\mu(A\sb{1})\leq \frac{1}{t}$. To ensure $\sum\sb{i=1}\sp{k}\lambda\sb{i}\chi\sb{A\sb{i}}\leq f$, there must be $\lambda\sb{1}\leq t-1$. Thus, it holds that
\[
\sum\sb{i=1}\sp{k}\lambda\sb{i}\mu(A\sb{i})= \lambda\sb{1}\mu(A\sb{1})\leq \frac{t-1}{t}<1.
\]
For the arbitrariness of the partition $\{A\sb{i}\}$, we conclude that $\int\sp{pan} fd\mu\leq 1$. In a similar way, we can prove that $\int\sp{pan} gd\mu\leq 1$ and $\int\sp{pan} (f+g)d\mu\leq 1$.
On the other hand, $(t-1)\chi\sb{\{t,t+1,t+2,\dots\}}\leq f$ and $(t-1)\mu(\{t,t+1,t+2,\dots\})=\frac{t-1}{t}$ for each $t\in X$. Thus we have that $\int\sp{pan} fd\mu\geq\sup\{\frac{t-1}{t}|t\in X\}=1$, and hence $\int\sp{pan} fd\mu=1$. Similarly, we have that $\int\sp{pan} gd\mu= 1$ and $\int\sp{pan} (f+g)d\mu= 1$.
\end{exa}

Notice that Ineq.(\ref{eq4-1}) holds for any monotone measure $\mu$. If $\mu$ is subadditive, then it becomes an equality. That is, the following result holds.

\begin{lemm}\label{lemm4-4}
Let $f, g\in {\cal F\sp{+}}$ with $\{f>0\}\cap~\{g>0\}=\emptyset$. If $\mu$ is subadditive, then
\begin{equation}\label{eq4-2}
\int\sp{pan}(f+g)d\mu= \int\sp{pan}fd\mu+\int\sp{pan}gd\mu.
\end{equation}
\end{lemm}

{\bf Proof.} By Lemma \ref{lemm4-1}, it suffices to prove that $\int\sp{pan}(f+g)d\mu\leq \int\sp{pan}fd\mu+\int\sp{pan}gd\mu$. For any given partition $\{A\sb{i}\}\sb{i=1}\sp{k}$ of $\{f+g>0\}$ and a sequence of positive number $\{\lambda\sb{i}\}\sb{i=1}\sp{k}$ such that $\sum\sb{i=1}\sp{k}\lambda\sb{i}\chi\sb{A\sb{i}}\leq f+g$.
Since $\{f>0\}\cap \{g>0\}=\emptyset$, we have that $\{A\sb{i}\cap \{f>0\}\}\sb{i=1}\sp{k}$(resp. $\{A\sb{i}\cap \{g>0\}\}\sb{i=1}\sp{k}$) is a partition of $\{f>0\}$(resp. $\{g>0\}$), and that
\[
\sum\sb{i=1}\sp{k}\lambda\sb{i}\chi\sb{A\sb{i}\cap \{f>0\}}\leq f,\qquad \sum\sb{i=1}\sp{k}\lambda\sb{i}\chi\sb{A\sb{i}\cap \{g>0\}}\leq g.
\]
Moreover, the subadditivity of $\mu$ implies that
\begin{eqnarray*}
\mu(A\sb{i})
&=&
\mu\left(A\sb{i}\bigcap(\{f>0\}\cup \{g>0\})\right)\\
&\leq&
\mu(A\sb{i}\cap \{f>0\})+ \mu(A\sb{i}\cap \{g>0\}).
\end{eqnarray*}
\begin{eqnarray*}
\sum\sb{i=1}\sp{k}\lambda\sb{i}\mu(A\sb{i})&\leq&\sum\sb{i=1}\sp{k}\lambda\sb{i}\mu(A\sb{i}\cap \{f>0\})+ \sum\sb{i=1}\sp{k}\lambda\sb{i}\mu(A\sb{i}\cap \{g>0\})\\
&\leq& \int\sp{pan}fd\mu+\int\sp{pan}gd\mu.
\end{eqnarray*}
Thus,
\begin{eqnarray*}
\int\sp{pan}(f+g) d\mu &=& \sup\{\sum\sb{i=1}\sp{n}\lambda\sb{i}\mu(A\sb{i}): \sum\sb{i=1}\sp{n}\lambda\sb{i}\chi\sb{A\sb{i}}\leq f+g,\\
 && \{A\sb{i}\}\sb{i=1}\sp{n}\subset {\mathcal F}
                  \mbox{ is a partition of } X, \lambda\sb{i}\geq 0,
                        n\in \mathbb{N}\}\\
                        &\leq& \int\sp{pan}fd\mu+\int\sp{pan}gd\mu.\Box
\end{eqnarray*}

\begin{prot}\label{prot4-5}
Let $(X, {\cal A}, \mu)$ be a subadditive monotone measure space and $f, g\in {\cal F\sp{+}}$. Then, $f+g$ is pan-integrable if and only if both $f$ and $g$ are pan-integrable.
\end{prot}

{\bf Proof.} Denote $A=\{x|f(x)\leq g(x)\}$ and $B=\{x|f(x)>g(x)\}$. Then $A\cap B=\emptyset$ and $f=f\cdot\chi\sb{A}+f\cdot\chi\sb{B}, g=g\cdot\chi\sb{A}+g\cdot\chi\sb{B}$. Moreover, $f\cdot\chi\sb{A}\leq g\cdot\chi\sb{A}$ and $g\cdot\chi\sb{B}\leq f\cdot\chi\sb{B}$. Thus,
\[
f+g\leq 2g\cdot\chi\sb{A}+2f\cdot\chi\sb{B}.
\]
Now, by the monotonicity, positive homogeneity and the positive sets disjointness additivity(Lemma \ref{lemm4-4})of the pan-integral, we conclude that
\begin{eqnarray*}
\int\sp{pan}(f+g)d\mu &\leq& \int\sp{pan}\left(2g\cdot\chi\sb{A}+2f\cdot\chi\sb{B}\right)d\mu\\
&=& \int\sp{pan}2g\cdot\chi\sb{A}d\mu+ \int\sp{pan}2f\cdot\chi\sb{B}d\mu\\
&=& 2\left(\int\sp{pan}g\cdot\chi\sb{A}d\mu+ \int\sp{pan}f\cdot\chi\sb{B}d\mu\right)\\
&\leq& 2\left(\int\sp{pan}g d\mu+ \int\sp{pan}f d\mu\right).
\end{eqnarray*}
Thus if $f+g$ is not pan-integrable then either $f$ or $g$ is not pan-integrable. One the other hand, if one of $f, g$ is not pan-integrable, then $f+g$ is not pan-integrable. Hence $f+g$ is pan-integrable if and only if both $f$ and $g$ are pan-integrable.$\Box$

\begin{nte}
By Proposition \ref{prot4-5}, if $f+g$ is not pan-integrable, then
\[
\max\left(\int\sp{pan}f d\mu, \int\sp{pan}g d\mu\right)=\infty.
\]
Thus, both sides of (\ref{eq3-1}) are equal to infinity and the equality holds.
\end{nte}

\section{The pan-integral for real-valued functions}

So far, the pan-integral is considered for $f\in {\cal M\sp{+}}$ only, whereas the Choquet integral for real-valued functions(not necessarily nonnegative) can be found in the literature. Let $(X, {\cal F}, \mu)$ be a given monotone measure space with $\mu(X)<\infty$. By $\bar{\mu}$ we denote the conjugation of $\mu$, i.e., $\bar{\mu}(A)=\mu(X)-\mu(X\setminus A)$ for each $A\in {\cal F}$. For any $f\in {\cal M}$, the asymmetric and the symmetric Choquet integrals are respectively defined by
 \[
\int\hspace{-0.4cm A}\hspace{0.4cm} f d\mu=\int f\sp{+} d\mu-\int f\sp{-} d\bar{\mu},\ \mbox{and}
 \]
 \[
\int\hspace{-0.4cm S}\hspace{0.4cm} f d\mu=\int f\sp{+} d\mu-\int f\sp{-} d\mu,
 \]
 where $f\sp{+}=\max(f, 0)$ and $f\sp{-}=\max(-f, 0)$. It is well-known that the asymmetric Choquet integral is comonotonic additivity and positive homogeneity, while the symmetric Choquet integral is homogneity but lacks comonotonic additivity \cite{Den94}. In contrast to the asymmetric Choquet integral, there is only one monotone measure in the definition of symmetric Choquet integral. Thus, to define the symmetric Choquet integral, we do not require the finiteness of $\mu$.

In the next section, we will investigate properties of the space of pan-integrable functions. In order to do this, here we introduce the pan-integral for real-valued(not necessarily positive) functions. Like the Choquet integral, there are two ways, namely, the so-called asymmetric one and the symmetric one. We will adopt the latter strategy.

\begin{dfnt}
Let $(X, {\cal A}, \mu)$ be a given monotone measure space, and $f\colon X\to \mathbb{R}$ a measurable function. If at least one of $\int\sb{X}\sp{pan}f\sp{+}d\mu, \int\sb{X}\sp{pan}f\sp{-}d\mu$ takes finite value, then we difine the pan-integral of $f$ (with respect to $\mu$) via
\begin{equation}\label{eq5-1}
\int\sb{X}\sp{pan}fd\mu=\int\sb{X}\sp{pan}f\sp{+}d\mu-\int\sb{X}\sp{pan}f\sp{-}d\mu.
\end{equation}
If both $\int\sb{X}\sp{pan}f\sp{+}d\mu$ and $\int\sb{X}\sp{pan}f\sp{-}d\mu$ are finite, then $\int\sb{X}\sp{pan}fd\mu<\infty$ and we say that $f$ is pan-integrable.

For $A\subset X$, define $\int\sb{A}\sp{pan}fd\mu=\int\sb{X}\sp{pan}f\cdot\chi\sb{A}d\mu$.
\end{dfnt}

\begin{exa}\label{exa5-2}
Let $X=\{x\sb{1}, x\sb{2}, x\sb{3}, x\sb{4}\}, {\cal A}={\cal P}(X)$ and the monotone measure $\mu$ be defined as follows:
\[\mu(\{x\sb{1}\})=\mu(\{x\sb{2}\})=1, \mu(\{x\sb{3}\})=2, \mu(\{x\sb{4}\})=\mu(\{x\sb{1}, x\sb{2}\})=1.5,\]
\[
 \mu(\{x\sb{1}, x\sb{3}\})=\mu(\{x\sb{2}, x\sb{4}\})=\mu(\{x\sb{3}, x\sb{4}\})=4, \mu(\{x\sb{1}, x\sb{4}\})=2.5,
\]
\[
\mu(\{x\sb{2}, x\sb{3}\})=3.5, \mu(\{x\sb{1}, x\sb{2}, x\sb{3}\})=\mu(\{x\sb{1}, x\sb{2}, x\sb{4}\})= 5,
\]
\[
\mu(\{x\sb{1}, x\sb{3}, x\sb{4}\})=4.5, \mu(\{x\sb{2}, x\sb{3}, x\sb{4}\})=6, \mu(X)=6.5.
\]
Let $f\colon X\to \mathbb{R}$ be defined as
\begin{eqnarray*}
 f(x)=\left \{
        \begin {array}{ll}
              2
              &\quad \mbox{if $x=x\sb{1}$,}\\
              -2
              &\quad \mbox{if $x=x\sb{2}$,}\\
              1
              &\quad \mbox{if $x=x\sb{3}$,}\\
              -1
              &\quad \mbox{if $x=x\sb{4}$,}\\
        \end {array}
       \right.
\end{eqnarray*}
Then $\int\sp{pan}f\sp{+}d\mu=\mu(\{x\sb{1},x\sb{3}\})=4$ and $\int\sp{pan}f\sp{-}d\mu=\mu(\{x\sb{2},x\sb{4}\})=4$. Thus $\int\sp{pan}fd\mu=\int\sp{pan}f\sp{+}d\mu-\int\sp{pan}f\sp{-}d\mu=0$.
\end{exa}

In the rest part of this section, we investigate some basic properties of the pan-integral.

\begin{prot}\label{prot5-3}
Let $(X, {\cal A})$ be a measurable space, $f, g\colon{\cal F}\to \mathbb{R}$ two ${\cal F}$-measurable functions and $c\in\mathbb{R}$ be a constant. Then we have the following:

(i) (Homogeneity) $\int\sp{pan} cf d\mu=c\int\sp{pan} f d\mu$;

(ii) (Monotonicity) $f\leq g$ implies that $\int\sp{pan} f d\mu\leq\int\sp{pan} f d\mu$;

(iii) If $|f|$ is integrable then $f$ is also integrable. Moreover, if $\mu$ is subadditive then $f$ is integrable if and only if $|f|$ is integrable.
\end{prot}

{\bf Proof.}\ The proofs of (i) and (ii) are standard and thus omitted. For (iii), notice first that $0\leq \max(f\sp{+}, f\sp{-})\leq |f|$. If $\int\sp{pan}|f| d\mu<\infty$, by the monotonicity of pan-integral, we then have that both $\int\sp{pan}f\sp{+} d\mu<\infty$ and $\int\sp{pan}f\sp{-} d\mu<\infty$, which imply the integrability of $f$.

By definition, the integrability of $f$ implies that both $f\sp{+}$ and $f\sp{-}$ are integrable. Suppose now $\mu$ is subadditive, Theorem \ref{th3-1} holds. Thus
\[
\int\sp{pan}|f| d\mu=\int\sp{pan}\left(f\sp{+}+ f\sp{-}\right) d\mu=\int\sp{pan}f\sp{+} d\mu+ \int\sp{pan}f\sp{-} d\mu<\infty. \Box
\]

The following result shows that the pan-integral is linear whenever the monotone measure $\mu$ is subadditive.

\begin{thrm}
Let $(X, {\cal A}, \mu)$ be a monotone measure space. If $\mu$ is subadditive, then the pan-integral with respect to $\mu$ is linear.
\end{thrm}

{\bf Proof.}\ By (i) of Proposition \ref{prot5-3}, it suffices to prove that $f+g$, the sum of two pan-integrable functions $f, g$, is pan-integrable, and
\[
\int\sp{pan}(f+g)d\mu= \int\sp{pan}fd\mu+\int\sp{pan}gd\mu.
\]
 We assume first that $\{f>0\}\cap \{g>0\}=\emptyset$. Then $(f+g)\sp{+}=f\sp{+}+ g\sp{+}$ and $(f+g)\sp{-}=f\sp{-}+ g\sp{-}$. By definition of the pan-integral and Lemma \ref{lemm4-4}, we have that
\begin{eqnarray*}
\int\sp{pan}(f+g)d\mu
&=&  \int\sp{pan}(f\sp{+}+g\sp{+})d\mu-\int\sp{pan}(f\sp{-}+g\sp{-})d\mu\\
&=& \left(\int\sp{pan}f\sp{+}d\mu+\int\sp{pan}g\sp{+}d\mu\right)-\left(\int\sp{pan}f\sp{-}d\mu+\int\sp{pan}g\sp{-}d\mu\right)\\
&=& \left(\int\sp{pan}f\sp{+}d\mu-\int\sp{pan}f\sp{-}d\mu\right)+\left(\int\sp{pan}g\sp{+}d\mu-\int\sp{pan}g\sp{-}d\mu\right)\\
&=& \int\sp{pan}fd\mu+ \int\sp{pan}gd\mu
\end{eqnarray*}
Now for two arbitrary pan-integrable functions $f, g$, denote
\[
A=\{x\in X| f(x)\geq 0, g(x)\geq 0\},
\]
\[
B=\{x\in X| f(x)< 0, g(x)< 0\},
\]
\[
C\sb{1}=\{x\in X| f(x)\geq 0, g(x)< 0, f(x)+g(x)\geq 0\},
\]
\[
C\sb{2}=\{x\in X| f(x)\geq 0, g(x)< 0, f(x)+g(x)< 0\},
\]
\[
D\sb{1}=\{x\in X| f(x)< 0, g(x)\geq 0, f(x)+g(x)\geq 0\},
\]
\[
D\sb{2}=\{x\in X| f(x)< 0, g(x)\geq 0, f(x)+g(x)< 0\}.
\]
Then $(f+g)\sp{+}=(f+g)\cdot\chi\sb{A}+(f+g)\cdot\chi\sb{C\sb{1}}+(f+g)\cdot\chi\sb{D\sb{1}}$. Observing that $A, B, C\sb{i}, D\sb{i}, i=1, 2$ are mutual disjoint sets, by the positive sets disjointness additivity, we have
\[
\int\sp{pan}(f+g)\sp{+}d\mu=\int\sp{pan}\left(f+g\right)\cdot\chi\sb{A}d\mu+\int\sp{pan}\left(f+g\right)\cdot\chi\sb{C\sb{1}}d\mu
+\int\sp{pan}\left(f+g\right)\cdot\chi\sb{D\sb{1}}d\mu.
\]
Similarly,
\[
\int\sp{pan}(f+g)\sp{-}d\mu=\int\sp{pan}-\left(f+g\right)\cdot\chi\sb{B}d\mu+\int\sp{pan}-\left(f+g\right)\cdot\chi\sb{C\sb{2}}d\mu
+\int\sp{pan}-\left(f+g\right)\cdot\chi\sb{D\sb{2}}d\mu.
\]
Noticing that both $f\cdot\chi\sb{A}$ and $g\cdot\chi\sb{A}$ are nonnegative, by Theorem \ref{th3-1} we have
\[
\int\sp{pan}\left(f+g\right)\cdot\chi\sb{A}d\mu=\int\sp{pan} f\cdot\chi\sb{A}d\mu+ \int\sp{pan} g\cdot\chi\sb{A}d\mu.
\]
As $0\leq \left(f+g\right)\cdot\chi\sb{C\sb{1}}\leq f\cdot\chi\sb{C\sb{1}}\leq f\sp{+}$ and $0\leq (-g)\cdot\chi\sb{C\sb{1}}\leq g\sp{-}$, they are pan-integrable. By using Theorem \ref{th3-1} one more time we get
\begin{eqnarray*}
\int\sp{pan}f\cdot\chi\sb{C\sb{1}}d\mu
&=&\int\sp{pan}\left(\left(f+g\right)\cdot\chi\sb{C\sb{1}}+(-g)\cdot\chi\sb{C\sb{1}}\right)d\mu\\
&=& \int\sp{pan}\left(f+g\right)\cdot\chi\sb{C\sb{1}} d\mu+\int\sp{pan}(-g)\cdot\chi\sb{C\sb{1}}d\mu.
\end{eqnarray*}
Thus
\[
\int\sp{pan}\left(f+g\right)\cdot\chi\sb{C\sb{1}} d\mu=\int\sp{pan}f\cdot\chi\sb{C\sb{1}}d\mu-\int\sp{pan}-g\cdot\chi\sb{C\sb{1}}d\mu.
\]
In a similar way, we can show the following equalities
\[
\int\sp{pan}\left(f+g\right)\cdot\chi\sb{D\sb{1}}d\mu=-\int\sp{pan}-f\cdot\chi\sb{D\sb{1}}d\mu+\int\sp{pan} g\cdot\chi\sb{D\sb{1}}d\mu,
\]
\[
\int\sp{pan}-\left(f+g\right)\cdot\chi\sb{B}d\mu=\int\sp{pan}-f\cdot\chi\sb{B}d\mu+ \int\sp{pan}-g\cdot\chi\sb{B}d\mu,
\]
\[
\int\sp{pan}-\left(f+g\right)\cdot\chi\sb{C\sb{2}}d\mu=-\int\sp{pan} f\cdot\chi\sb{C\sb{2}} d\mu+ \int\sp{pan} -g\cdot\chi\sb{C\sb{2}} d\mu,
\]
\[
\int\sp{pan}-\left(f+g\right)\cdot\chi\sb{D\sb{2}}d\mu=\int\sp{pan} -f\cdot\chi\sb{D\sb{2}} d\mu- \int\sp{pan} g\cdot\chi\sb{D\sb{2}} d\mu.
\]
Combining these equalities, by (i) of Proposition \ref{prot5-3} and the positive sets disjointness additivity, we can get the desired result.$\Box$

\begin{remk}
We have proved that the pan-integral for nonnegative functions possesses the positive sets disjointness superadditivity (see Lemma \ref{lemm4-1}). For general measurable functions(not necessarily nonnegative), however, this property fails. Let us reconsider Example \ref{exa5-2}. Let $g=f\cdot\chi\sb{\{x\sb{1},x\sb{4}\}}$ and $h=f\cdot\chi\sb{\{x\sb{2},x\sb{3}\}}$. Then $\{g>0\}\cap \{h>0\}=\emptyset$. It is easy to see that $\int\sp{pan}g d\mu=2\mu(\{x\sb{1}\})-\mu(\{x\sb{4}\})=0.5, \int\sp{pan} h d\mu=\mu(\{x\sb{3}\})-2\mu(\{x\sb{2}\})=0$. That is,
\[
\int\sp{pan}(g+h) d\mu=\int\sp{pan}f d\mu=0<0.5=\int\sp{pan}g d\mu+ \int\sp{pan}h d\mu.
\]
\end{remk}

We close this section by presenting some convergence theorem for pan-integral. The proofs of these results are similar to that of the Lebesgue integral and thus are omitted.

\begin{prot}\label{prot5-6}
(Levi) Let $f\sb{k}\colon E\to [0, \infty]$ be a sequence of nonnegative measurable functions such that $f\sb{k}\leq f\sb{k+1}$ for each $k$. Let $\lim\sb{k\to\infty}f\sb{k}(x)=f(x), x\in E$. Then
\[
\lim\sb{k\to\infty}\int\sp{pan}\sb{E} f\sb{k} d\mu= \int\sp{pan}\sb{E} f d\mu.
\]
\end{prot}

By Proposition \ref{prot5-6}, the following result can be proved.

\begin{prot}\label{prot5-7}
(Fatou)  Let $f\sb{k}\colon E\to [0, \infty]$ be a sequence of nonnegative measurable functions. Then
\[
\int\sp{pan}\sb{E}\liminf\sb{k\to\infty} f\sb{k} d\mu\leq \liminf\sb{k\to\infty}\int\sp{pan}\sb{E} f\sb{k} d\mu.
\]
\end{prot}

\section{Pan-integrable functional space}
Let $(X, {\cal A}, \mu)$ be a monotone measure space and $E\subset X$ a measurable set.
For $p\geq 1$, denote
\[
{\cal L}\sb{\mu}\sp{p}(E)=\left\{f\colon E\to\mathbb{R}\colon\left(\int\sp{pan}\sb{E}|f|\sp{p} d\mu\right)\sp{\frac{1}{p}}<\infty\right\}.
\]
By using Proposition \ref{prot5-3}(iii) and Theorem \ref{th3-1}, it can be shown that
\begin{thrm}
Let $(X, {\cal A}, \mu)$ be a monotone measure space. If $\mu$ is subadditive, then for any $f, g\in{\cal L}\sb{\mu}\sp{p}(E)$ and any real numbers $\alpha, \beta$, we have that $\alpha f+\beta g\in{\cal L}\sb{\mu}\sp{p}(E)$, i.e., ${\cal L}\sb{\mu}\sp{p}(E)$ is a linear space.
\end{thrm}

We will prove that ${\cal L}\sb{\mu}\sp{p}(E)$ is a Banach space if $\mu$ is \emph{subadditive and continuous from below}.
Due to Theorem \ref{th3-1}, most of the results can be proved in a similar way of the Lebesgue integral. So, we will omit the arguments. To prove ${\cal L}\sb{\mu}\sp{p}(E)$ is a complete metric space, we need the following concept.

Let $f\colon E\to [0, \infty]$ be a measurable function and $p\geq 1$. Define
\[
\|f\|\sb{\mu,p}=\left(\int\sp{pan}\sb{E}|f|\sp{p} d\mu\right)\sp{\frac{1}{p}},
\]
and
\[
\rho(f, g)= \|f-g\|\sb{\mu,p}.
\]
Then $\rho(f, g)$ is a metric whenever $\mu$ is subadditive and continuous from below. In fact

(i) for any monotone measure, $\rho(f, g)=\rho(g, f)$;

(ii) if $\mu$ is continuous from below then by (i) of Note \ref{nte2-4} we know that $\rho(f, g)=0$ if and only if $f=g$~ $\mu$-a.e. on $E$;

(iii) if $\mu$ is subadditive and continuous from below, then the triangle inequality follows from the Minkowski inequality below.

\begin{lemm}\label{lem6-2}(see also \cite{YanOuysubmitted})
Let $(X, {\cal A}, \mu)$ be a subadditive monotone measure space.

(i)(H\"{o}lder) If $p, q\geq 1$ satisfy $\frac{1}{p}+\frac{1}{q}=1$ and $\|f\|\sb{\mu,p}<\infty, \|g\|\sb{\mu,q}<\infty$, then
\[
\|fg\|\sb{\mu,1}\leq\|f\|\sb{\mu,p}\|g\|\sb{\mu,q}.
\]

(ii)(Minkowski) For $p\geq 1$ and $f, g\in{\cal L}\sb{\mu}\sp{p}(E)$, we have that
\[
\|f+g\|\sb{\mu,p}\leq\|f\|\sb{\mu,p}+\|g\|\sb{\mu,p}.
\]
\end{lemm}

By Proposition \ref{prot5-7} and \ref{lem6-2}, we can prove the completeness of ${\cal L}\sb{\mu}\sp{p}(E)$

\begin{thrm}
Let $(X, {\cal A}, \mu)$ be a monotone measure space and $E\in{\cal F}$. If $\mu$ is subadditive and continuous from below, then the set ${\cal L}\sb{\mu}\sp{p}(E)$ is a complete metric space(equipped with the metric $\rho(f, g)$).
\end{thrm}

\section{Conclusions}
In this paper, we have investigated the properties of $+,\cdot$-based pan-integral. It has been proved that this integral is linear whenever the monotone measure is subadditive. As a consequence, many results of the Lebesgue integral theory including the $L\sp{p}$ space, can be generalized to the pan-integral. Since the concave integral \cite{LehTep08,Leh09,Tep09} coincides with the pan-integral whenever $\mu$ is subadditive \cite{MesLiOuysubmitted,ZhaYanOuy16}, the results obtained in this paper also hold for concave integral.

It should be pointed out that an outer measure is subadditive. Thus we can define a Lebesgue-like integral(possesses linearity) from an outer measure, and the $L\sp{p}$ theory holds for this integral. Moreover, the domain of an outer measure is the power set, measurability restrictions of functions are thus unnecessary for this integral.

In future work, we will investigate similar results for $(\oplus,\otimes)$-based pan-integral. We stress that this is not a trivial generalization, since we even do not know whether or not the positive homogeneity holds for general pan-integrals. On the other hand, the subadditivity is a rather strong restriction for monotone measure. It is of interest to exam the results of this paper under some weaker conditions, such as the uniform autocontinuity \cite{WanKli09}.

\section*{Acknowledgment}
This work was partially supported by NSFC (nos.11371332 and 11571106) and the NSF of Zhejiang Province(no.LY15A010013), and by the grant APVV-14-0013.

\end{document}